\documentclass[10pt]{bmc_article}
\usepackage{amsmath,amssymb,amscd,enumerate,latexsym,graphicx}
\usepackage{cite} 
\usepackage{url}  
\usepackage{ifthen}  
\usepackage{multicol}   
\usepackage[utf8]{inputenc} 
\usepackage{amsfonts,epsfig}
\usepackage{pifont,indentfirst,subfigure}
\usepackage{amssymb,amsmath,amsthm,graphicx,amsxtra}

\newtheorem{thm}{Theorem}
\newtheorem{lem}{Lemma}
\newtheorem{rmk}{Remark}

\newtheorem{example}{Example}

\def\includegraphics{}

\newboolean{publ}

\newenvironment{bmcformat}{\baselineskip20pt\sloppy\setboolean{publ}{false}}{\baselineskip20pt\sloppy}

\begin{document}
\begin{bmcformat}
\title{Lower bounds of the minimum eigenvalue for $M$-matrices}
\author{Jianxing Zhao\correspondingauthor%
        \email{Jianxing Zhao\correspondingauthor - zjx810204@163.com}
         and
        Caili Sang%
         \email{Caili Sang - sangcl@126.com}
}
\address{%
College of Science, Guizhou Minzu University, Guiyang, Guizhou
550025, P.R. China
}%

\maketitle
\begin{abstract}
Some monotone increasing sequences of the lower bounds for the minimum eigenvalue of $M$-matrices are given. It is proved that these sequences are convergent and improve some existing results. Numerical examples show that these sequences are more accurate than some existing results and could reach the true value of the minimum eigenvalue in some cases.\\
\textbf {MSC:} 15A06; 15A15; 15A48\\
\textbf {Keywords:} $M$-matrix; nonnegative matrix; Hadamard product; spectral radius; minimum
eigenvalue

\end{abstract}

\ifthenelse{\boolean{publ}}{\begin{multicols}{2}}{}
\section{Introduction}

For a positive integer $n~(n\geq 2)$, $N$ denotes the set $\{1, 2, \ldots, n\}$,
and $\mathbb{R}^{n\times n}(\mathbb{C}^{n\times n})$ denotes the set
of all ${n\times n}$ real (complex) matrices throughout.
For $A=[a_{ij}]\in \mathbb{R}^{n\times n}$, we write $A\geq 0$ if $a_{ij}\geq 0,i,j\in N.$
If $A\geq 0$, we say $A$ is nonnegative.

A matrix $A=[a_{ij}]\in \mathbb{R}^{n\times n}$ is called a nonsingular
$M$-matrix if $a_{ij}\leq 0, i\neq j,i,j\in N$ and the inverse of $A$, denoted by $A^{-1}$, is nonnegative. Denote by $M_n$ the set of all
$n\times n$ nonsingular $M$-matrices (see \cite{Berman}). If $A$ is a nonsingular $M$-matrix, then there exists a positive
eigenvalue of $A$ equal to $\tau(A)=\rho(A^{-1})^{-1}$, where
$ \rho(A^{-1})$ is the perron eigenvalue of the nonnegative matrix
$A^{-1}.$  It is easy to prove that
$\tau(A)=\min\{|\lambda|:\lambda\in\sigma(A)\},$ where $\sigma(A)$
denotes the spectrum of  $A$. $\tau(A)$ is called the minimum eigenvalue of $A$ (see \cite{Horn}).
If $G$ is the diagonal matrix of an $M$-matrix $A$,
then the spectral radius of the Jacobi iterative matrix $J_A=G^{-1}(G-A)$ of $A$,
denoted by $\rho(J_A)$, is less than 1 (see \cite{Berman}).

For two real matrices $A=[a_{ij}]$ and $B=[b_{ij}]$ of the same size, the Hadamard product of $A$ and $B$ is defined as the
matrix $A\circ B=[a_{ij}b_{ij}]$. If $A\in M_n$ and $B \geq 0$, then it is clear
that $B\circ A^{-1}\geq 0$ (see \cite{Horn}).

Let $A=[a_{ij}]\in \mathbb{R}^{n\times n},a_{ii}\neq 0,i\in N$, and $A^{-1}=[\alpha_{ij}]$. For $i,j,k\in N,j\neq i,t=1,2,\ldots,$ denote
\[d_i=\frac{\sum\limits_{j\neq i}{|a_{ij}|}}{|a_{ii}|},d=\max\limits_{i\in N}d_i,\varphi_i=\frac{1}{a_{ii}-\sum\limits_{k\neq i}|a_{ik}|d_k};
r_i=\max\limits_{j\neq{i}}\bigg\{\frac{|a_{ji}|}{|a_{jj}|-\sum\limits_{k \neq j,
i}|a_{jk}|}\bigg\}, m_{ji}=\frac{|a_{ji}|+\sum\limits_{k\neq{j,i}}|a_{jk}|r_i}{|a_{jj}|},\]
\[h_i=\max\limits_{j\neq{i}}\bigg\{\frac{|a_{ji}|}{|a_{jj}|m_{ji}-\sum\limits_{k\neq{j,i}}|a_{jk}|m_{ki}}\bigg\},
u_{ji}=\frac{|a_{ji}|+\sum\limits_{k\neq{j,i}}|a_{jk}|m_{ki}h_i}{|a_{jj}|},
u_i=\max\limits_{j\neq{i}}\{u_{ij}\}.\]
\[
u^{(0)}_{ji}=u_{ji},~p_{ji}^{(t)}=\frac{|a_{ji}|+\sum\limits_{k\neq{j,i}}|a_{jk}|u_{ki}^{(t-1)}}{|a_{jj}|},
~p^{(t)}_i=\max\limits_{j\neq{i}}\{p^{(t)}_{ij}\},\]
\[h^{(t)}_i=\max\limits_{j\neq{i}}\bigg\{\frac{|a_{ji}|}{|a_{jj}|p^{(t)}_{ji}-\sum\limits_{k\neq{j,i}}|a_{jk}|p^{(t)}_{ki}}\bigg\},
~u^{(t)}_{ji}=\frac{|a_{ji}|+\sum\limits_{k\neq{j,i}}|a_{jk}|p^{(t)}_{ki}h^{(t)}_i}{|a_{jj}|};
\phi_i^{(t)}=\frac{1}{a_{ii}-\sum\limits_{j\neq i}|a_{ij}|p_{ji}^{(t)}}.
\]

Recall that $A=[a_{ij}]\in\mathbb{C}^{n\times n}$ is called diagonally dominant if $d_i\leq 1$ for all $i\in N$. If $d_i<1$, we say that $A$ is strictly diagonally dominant. It is well known that a strictly diagonally dominant matrix is nonsingular. $A$ is called weakly chained
diagonally dominant if $d_{i}\leq 1, J(A)=\{i\in N: d_i<1 \}\neq \varnothing$ and for all $i\in N/J(A),$
there exist indices $i_1,i_2,\ldots,i_k$ in $N$ with $a_{i_li_{l+1}}\neq 0,0\leq l\leq k-1$, where $i_0=i$ and
$i_k\in J(A)$. Notice that a strictly diagonally dominant matrix is also weakly chained
diagonally dominant (see \cite{Shi}).

Estimating the bounds for the minimum eigenvalue of $M$-matrices is an interesting
subject in matrix theory, it has important applications in many practical problems
(see [3-12]) and various refined bounds can be found in [3-8]. Hence, it is necessary to estimate the bounds for $\tau(A)$.

In \cite{Shi}, Shivakumar \emph{et al}. gave the following bounds for $\tau(A)$: Let $A=[a_{ij}]\in M_n$ be
weakly chained diagonally dominant and $A^{-1}=[\alpha_{ij}]$. Then
\begin{eqnarray}\label{shi-equ}
\min\limits_{i\in N}\sum\limits_{j=1}^n{a_{ij}}\leq\tau(A)\leq \max\limits_{i\in N}\sum\limits_{j=1}^n{a_{ij}},\tau(A)\leq\min\limits_{i\in N}{a_{ii}}~\textmd{and}~\frac{1}{\max\limits_{i\in N}\sum\limits_{j=1}^n{\alpha_{ij}}}\leq \tau(A)\leq\frac{1}{\min\limits_{i\in N}\sum\limits_{j=1}^n{\alpha_{ij}}}.
\end{eqnarray}

Subsequently, Tian and Huang \cite{tian} obtained a lower bound for $\tau(A)$ using the spectral
radius of the Jacobi iterative matrix $J_A$ of $A$: Let $A=[a_{ij}]\in M_n$ and $A^{-1}=[\alpha_{ij}]$. Then
\begin{eqnarray}\label{tgx-equ1}
\tau(A)\geq \frac{1}{[1+(n-1)\rho(J_A)]\max\limits_{i\in N}{\alpha_{ii}}}.
\end{eqnarray}
Furthermore, when $A$ is a strictly diagonally dominant $M$-matrix, they provided
lower bound for $\tau(A)$ which depend only on the entries of $A$: If $A=[a_{ij}]\in M_n$ is strictly diagonally dominant, then
\begin{eqnarray}\label{tgx-equ2}
\tau(A)\geq \frac{1}{[1+(n-1)d]\max\limits_{i\in N}{\varphi_{i}}}.
\end{eqnarray}

In 2013, Li \emph{et al}. \cite{lcq} improved (\ref{tgx-equ1}) and (\ref{tgx-equ2}), and presented the following result: Let $A=[a_{ij}]\in M_n$ and $A^{-1}=[\alpha_{ij}]$. Then
\begin{eqnarray}\label{lcq-equ1}
\tau(A)\geq\frac{2}{\max\limits_{i\neq j}\big\{\alpha_{ii}+\alpha_{jj}+[(\alpha_{ii}-\alpha_{jj})^2+4(n-1)^2\alpha_{ii}\alpha_{jj}\rho^2(J_A)]^{\frac{1}{2}}\big\}}.
\end{eqnarray}
Furthermore, when $A$ is a strictly diagonally dominant $M$-matrix, they also obtained
lower bound for $\tau(A)$ which depend only on the entries of $A$: If $A=[a_{ij}]\in M_n$ is strictly diagonally dominant, then
\begin{eqnarray}\label{lcq-equ2}
\tau(A)\geq\frac{2}{\max\limits_{i\neq j}\big\{\varphi_{i}+\varphi_{j}+[\varphi_{ij}^2+4(n-1)^2\varphi_{i}\varphi_{j}d^2]^{\frac{1}{2}}\big\}},
\end{eqnarray}
where
$\varphi_{ij}=\max\{\varphi_i,\varphi_j\}-\min\{a_{ii}^{-1},a_{jj}^{-1}\}.$

In 2015, Wang and Sun \cite{wf} gave the following result:
Let $A=[a_{ij}]\in M_n$ and $A^{-1}=[\alpha_{ij}]$. Then
\begin{eqnarray}\label{wf-equ1}
\tau(A)\geq\frac{2}{\max\limits_{i\neq j}\big\{\alpha_{ii}+\alpha_{jj}+[(\alpha_{ii}-\alpha_{jj})^2+4(n-1)^2\alpha_{ii}\alpha_{jj}u_iu_j]^{\frac{1}{2}}\big\}}.
\end{eqnarray}

Recently, Zhao and Sang \cite{zjx2016mineigenvaule} obtained the following result:
Let $A=[a_{ij}]\in M_n$ and $A^{-1}=[\alpha_{ij}]$.
Then, for $t=1,2,\ldots$,
\begin{eqnarray}\label{zjx2016-equ1}
 \tau(A)\geq \frac{2}{\max\limits_{i\neq j}\Big\{\alpha_{ii}+\alpha_{jj}+\Big[(\alpha_{ii}-\alpha_{jj})^2+4(n-1)^2p_i^{(t)}p_j^{(t)}\alpha_{ii}\alpha_{jj}
 \Big]^{\frac{1}{2}}\Big\}}=\Upsilon_t.
\end{eqnarray}
Similarly, they presented lower bounds for $\tau(A)$ which depend only on the entries of $A$ in the case of
$A$ is a strictly diagonally dominant $M$-matrix:
If $A=[a_{ij}]\in M_n$ is strictly diagonally dominant, then for $t=1,2,\ldots$,
\begin{eqnarray}\label{zjxcor1-equ1}
 \tau(A)\geq \frac{2}{\max\limits_{i\neq j}\Big\{\phi_i^{(t)}+\phi_j^{(t)}+\Big[(\psi_{ij}^{(t)})^2+4(n-1)^2p_i^{(t)}p_j^{(t)}\phi_{i}^{(t)}\phi_{j}^{(t)}
 \Big]^{\frac{1}{2}}\Big\}}=\widetilde{\Upsilon}_t,
\end{eqnarray}
where $\psi_{ij}^{(t)}=\max\{\phi_i^{(t)},\phi_j^{(t)}\}-\min\{a_{ii}^{-1},a_{jj}^{-1}\}.$

Next, we continue to research the problems mentioned above and give several convergent sequences of the lower bounds for $\tau(A)$.
Numerical examples show that the new lower bounds are more accurate than these lower bounds obtained by inequalities (1)-(8).

\section{Some lemmas}
In this section, we give some lemmas, which will be useful in the following proofs.

\begin{lem}\label{hadamardxy}\emph{\cite{Horn}}
Let $A,B\in \mathbb{R}^{n\times n}$, and let $X,Y\in\mathbb{R}^{n\times n}$ be diagonal matrices. Then
\[X(A\circ B)Y=(XAY)\circ B =(XA)\circ (BY)=(AY)\circ(XB)=A\circ(XBY).\]
\end{lem}

\begin{lem}\label{byp} \emph{\cite{Horn}}
Let  $A=[a_{ij}]\in \mathbb{C}^{n\times n}$. Then all the eigenvalues of  $A$ lie
in the region
$$\bigcup\limits_{i,j\in N,i\neq j}\Big\{z\in \mathbb{C} :|z-a_{ii}||z-a_{jj}| \leq \sum\limits_{k \neq  i}
|a_{ki}|\sum\limits_{k \neq j}|a_{kj}|\Big\}.$$
\end{lem}

\begin{lem}\label{zjx2015-1lemma1}\emph{\cite{zjx2016mineigenvaule}}
If $A=[a_{ij}]\in M_n$ is strictly diagonally dominant, then
$A^{-1}=[\alpha_{ij}]$ exists, and for all
$i,j\in{N},j\neq{i},t=1,2,\ldots,$
\begin{eqnarray*}
&&\emph{(a)}~1>r_i\geq m_{ji}\geq u_{ji}=u^{(0)}_{ji}\geq{p^{(1)}_{ji}}\geq {u^{(1)}_{ji}}\geq{p^{(2)}_{ji}}\geq{u^{(2)}_{ji}}\geq\ldots\geq{p^{(t)}_{ji}}\geq{u^{(t)}_{ji}}\geq\ldots\geq0;\hspace*{6cm}\\
&&\emph{(b)}~\alpha_{ji}\leq p_{ji}^{(t)}\alpha_{ii}; ~\frac{1}{a_{ii}}\leq \alpha_{ii}\leq \phi_i^{(t)}.
\end{eqnarray*}
\end{lem}

\begin{lem}\label{hzg}\emph{\cite{hzg}}
If $A=[a_{ij}]\in M_n$ is strictly diagonally dominant, then
$A^{-1}=[\alpha_{ij}]$ exists, and for all
$i\in{N},$
$\alpha_{ii}\geq\frac{1}{a_{ii}-\sum\limits_{k\neq i}\frac{a_{ik}a_{ki}}{a_{kk}}}.$
\end{lem}

\begin{lem}\label{ssj}\emph{\cite{cfb2009}}
If $A^{-1}$ is a doubly stochastic matrix, then $Ae=e, A^Te=e,$ where $e=[1,1,\ldots,1]^T.$
\end{lem}

\section{Main results}
In this section, we present our main results.
\begin{thm}\label{th1}
Let $A=[a_{ij}]\in M_n, B=[b_{ij}]\geq 0$, and $A^{-1}=[\alpha_{ij}]$.
Then, for $t=1,2,\ldots$,
\begin{eqnarray*}
 \rho(B\circ A^{-1})&\leq&
 \max\limits_{i\neq j}\frac{1}{2}\Big\{b_{ii}\alpha_{ii}+b_{jj}\alpha_{jj}+\Big[(b_{ii}\alpha_{ii}-b_{jj}\alpha_{jj})^2
 +4\alpha_{ii}\alpha_{jj}\Big(\sum\limits_{k\neq i}b_{ki}p_{ki}^{(t)}\Big)\Big(\sum\limits_{k\neq
j}b_{kj}p_{kj}^{(t)}\Big)\Big]^{\frac{1}{2}}\Big\} \\
 &\leq&\max\limits_{i\in N}\Big\{\big(b_{ii}+\sum\limits_{k\neq i}b_{ki}p_{ki}^{(t)}\big)\alpha_{ii}\Big\}.
 \end{eqnarray*}
\end{thm}

\begin{proof} (a) Since $A$ is an $M$-matrix, there exists a positive diagonal matrix $X$, such that $X^{-1}AX$ is a strictly diagonally dominant $M$-matrix (see \cite{Horn}), and, by Lemma \ref{hadamardxy},  
$$\rho(B\circ A^{-1})=\rho(X^{-1}(B\circ A^{-1})X)=\rho(B\circ(X^{-1}AX)^{-1}).$$
Hence, for convenience and without loss of generality, we assume that $A$ is a strictly diagonally
dominant matrix.

Let $\lambda=\rho(B\circ A^{-1})$, then $\lambda \geq b_{ii}\alpha_{ii},\forall i\in N.$ By
Lemma \ref{byp} and Lemma \ref{zjx2015-1lemma1}, there are $i,j\in N,~i\neq j$ such that
\begin{eqnarray*}
|\lambda-b_{ii}\alpha_{ii}||\lambda-b_{jj}\alpha_{jj}|&\leq&
\Big(\sum\limits_{k\neq i}b_{ki}\alpha_{ki}\Big)
\Big(\sum\limits_{k\neq
j}b_{kj}\alpha_{kj}\Big)
\leq\Big(\sum\limits_{k\neq i}b_{ki}p_{ki}^{(t)}\alpha_{ii}\Big)\Big(\sum\limits_{k\neq
j}b_{kj}p_{kj}^{(t)}\alpha_{jj}\Big)\\
&=&\alpha_{ii}\alpha_{jj}\Big(\sum\limits_{k\neq i}b_{ki}p_{ki}^{(t)}\Big)\Big(\sum\limits_{k\neq
j}b_{kj}p_{kj}^{(t)}\Big),
\end{eqnarray*}
i.e.,
\begin{eqnarray}\label{th1-equ1}
(\lambda-b_{ii}\alpha_{ii})(\lambda-b_{jj}\alpha_{jj})\leq
\alpha_{ii}\alpha_{jj}\Big(\sum\limits_{k\neq i}b_{ki}p_{ki}^{(t)}\Big)\Big(\sum\limits_{k\neq
j}b_{kj}p_{kj}^{(t)}\Big).
\end{eqnarray}
From (\ref{th1-equ1}), we have
\begin{eqnarray}\nonumber
 \lambda\leq \frac{1}{2}\Big\{b_{ii}\alpha_{ii}+b_{jj}\alpha_{jj}+\Big[(b_{ii}\alpha_{ii}-b_{jj}\alpha_{jj})^2
 +4\alpha_{ii}\alpha_{jj}\Big(\sum\limits_{k\neq i}b_{ki}p_{ki}^{(t)}\Big)\Big(\sum\limits_{k\neq
j}b_{kj}p_{kj}^{(t)}\Big)\Big]^{\frac{1}{2}}\Big\},
\end{eqnarray}
that is,
\begin{eqnarray*}
 \rho(B\circ A^{-1})&\leq& \frac{1}{2}\Big\{b_{ii}\alpha_{ii}+b_{jj}\alpha_{jj}+\Big[(b_{ii}\alpha_{ii}-b_{jj}\alpha_{jj})^2
 +4\alpha_{ii}\alpha_{jj}\Big(\sum\limits_{k\neq i}b_{ki}p_{ki}^{(t)}\Big)\Big(\sum\limits_{k\neq
j}b_{kj}p_{kj}^{(t)}\Big)\Big]^{\frac{1}{2}}\Big\}\\
&\leq&\max\limits_{i\neq j}\frac{1}{2}\Big\{b_{ii}\alpha_{ii}+b_{jj}\alpha_{jj}+\Big[(b_{ii}\alpha_{ii}-b_{jj}\alpha_{jj})^2
 +4\alpha_{ii}\alpha_{jj}\Big(\sum\limits_{k\neq i}b_{ki}p_{ki}^{(t)}\Big)\Big(\sum\limits_{k\neq
j}b_{kj}p_{kj}^{(t)}\Big)\Big]^{\frac{1}{2}}\Big\}.
\end{eqnarray*}

(b) Without loss of generality, for $i,j\in N,i\neq j,$ assume that
\[b_{jj}\alpha_{jj}+\alpha_{jj}\sum\limits_{k\neq j}b_{kj}p_{kj}^{(t)}\leq b_{ii}\alpha_{ii}+\alpha_{ii}\sum\limits_{k\neq i}b_{ki}p_{ki}^{(t)},\]
i.e.,
\[
\alpha_{jj}\sum\limits_{k\neq j}b_{kj}p_{kj}^{(t)}\leq b_{ii}\alpha_{ii}-b_{jj}\alpha_{jj}+\alpha_{ii}\sum\limits_{k\neq i}b_{ki}p_{ki}^{(t)}.
\]
Let $\Delta_{ij}=\Big[(b_{ii}\alpha_{ii}-b_{jj}\alpha_{jj})^2+4\alpha_{ii}\alpha_{jj}\Big(\sum\limits_{k\neq i}b_{ki}p_{ki}^{(t)}\Big)\Big(\sum\limits_{k\neq j}b_{kj}p_{kj}^{(t)}\Big)\Big]^{\frac{1}{2}}$. Then
\begin{eqnarray*}
\Delta_{ij}&\leq& \Big[(b_{ii}\alpha_{ii}-b_{jj}\alpha_{jj})^2+4\alpha_{ii}\Big(\sum\limits_{k\neq i}b_{ki}p_{ki}^{(t)}\Big)\Big(b_{ii}\alpha_{ii}-b_{jj}\alpha_{jj}+\alpha_{ii}\sum\limits_{k\neq i}b_{ki}p_{ki}^{(t)}\Big)\Big]^{\frac{1}{2}}\\
 &=&\Big[\Big(b_{ii}\alpha_{ii}-b_{jj}\alpha_{jj}+2\alpha_{ii}\sum\limits_{k\neq i}b_{ki}p_{ki}^{(t)}\Big)^2\Big]^{\frac{1}{2}}\\
 &=&b_{ii}\alpha_{ii}-b_{jj}\alpha_{jj}+2\alpha_{ii}\sum\limits_{k\neq i}b_{ki}p_{ki}^{(t)}.
\end{eqnarray*}
Further, we have
\begin{eqnarray*}
b_{ii}\alpha_{ii}+b_{jj}\alpha_{jj}+\Delta_{ij}\leq 2b_{ii}\alpha_{ii}+2\alpha_{ii}\sum\limits_{k\neq i}b_{ki}p_{ki}^{(t)},
\end{eqnarray*}
then
\[
 \rho(B\circ A^{-1})\leq\max\limits_{i\neq j}\frac{1}{2}\{b_{ii}\alpha_{ii}+b_{jj}\alpha_{jj}+\Delta_{ij}\}
 \leq\max\limits_{i\in N}\Big\{\Big(b_{ii}+\sum\limits_{k\neq i}b_{ki}p_{ki}^{(t)}\Big)\alpha_{ii}\Big\}.
\]
The proof is completed.
\end{proof}

\begin{thm}\label{th2} Let $A=[a_{ij}]\in M_n$ and $A^{-1}=[\alpha_{ij}]$.
Then, for $t=1,2,\ldots$,
\begin{eqnarray}\label{th2-equ1}
 \tau(A)\geq\frac{2}{\max\limits_{i\neq j}\Big\{\alpha_{ii}+\alpha_{jj}+\Big[(\alpha_{ii}-\alpha_{jj})^2+4\alpha_{ii}\alpha_{jj}\sum\limits_{k\neq i}p_{ki}^{(t)}\sum\limits_{k\neq
j}p_{kj}^{(t)}
 \Big]^{\frac{1}{2}}\Big\}}=\Gamma_t.
\end{eqnarray}
\end{thm}
\begin{proof}
Let all entries of $B$ in Theorem \ref{th1} be 1. Then
\begin{eqnarray}\label{th2-equ3}
\rho(A^{-1})\leq\max\limits_{i\neq j}\frac{1}{2}\Big\{\alpha_{ii}+\alpha_{jj}+\Big[(\alpha_{ii}-\alpha_{jj})^2
 +4\alpha_{ii}\alpha_{jj}\sum\limits_{k\neq i}p_{ki}^{(t)}\sum\limits_{k\neq
j}p_{kj}^{(t)}\Big]^{\frac{1}{2}}\Big\}.
\end{eqnarray}
From inequality (\ref{th2-equ3}) and $\tau(A)=\frac{1}{\rho(A^{-1})}$, the conclusion follows obviously.
\end{proof}

Similar to the proof of Theorem \ref{th2}, the following theorem is obtained easily.
\begin{thm}\label{th3} Let $A=[a_{ij}]\in M_n$ and $A^{-1}=[\alpha_{ij}]$.
Then, for $t=1,2,\ldots$,
\begin{eqnarray*}
 \tau(A)\geq \frac{1}{\max\limits_{i\in N}\Big\{\big(1+\sum\limits_{k\neq i}p_{ki}^{(t)}\big)\alpha_{ii}\Big\}}=\Omega_t.
\end{eqnarray*}
\end{thm}

\begin{thm}\label{th4}
The sequence $\{\Gamma_t\}$ $(\{\Omega_t\}),t=1,2,\ldots$ obtained
from Theorem \ref{th2} (Theorem \ref{th3}) is monotone increasing with an upper bound
$\tau(A)$ and, consequently, is convergent.
\end{thm}
\begin{proof}
By Lemma \ref{zjx2015-1lemma1}, we have $1>p^{(t)}_{ji}\geq p^{(t+1)}_{ji}\geq
0,j,i\in N,j \neq i, t=1,2,\ldots$. Thus, $\{\Gamma_t\}$ ($\{\Omega_t\}$) is monotonically
increasing sequence. Hence, the
sequence $\{\Gamma_t\}$ ($\{\Omega_t\}$) is convergent.
\end{proof}

\begin{rmk}\label{rmk1}
From Theorem \ref{th1} and the proof of Theorem \ref{th2}, it is easily to see that
if $A=[a_{ij}]\in M_n$ and $A^{-1}=[\alpha_{ij}]$, then $\tau(A)\geq\Gamma_t\geq \Omega_t, t=1,2,\ldots.$
\end{rmk}

Let $A$ is a strictly diagonally dominant $M$-matrix. Then two new lower bounds for $\tau(A)$, which depend only on the entries of $A$, are obtained .
\begin{thm}\label{th5}
If $A=[a_{ij}]\in M_n$ is strictly diagonally dominant, then for $t=1,2,\ldots$,
\begin{eqnarray}\label{th5-equ1}
 \tau(A)\geq \frac{2}{\max\limits_{i\neq j}\Big\{\phi_i^{(t)}+\phi_j^{(t)}+\Big[(\phi_{ij}^{(t)})^2+4\phi_{i}^{(t)}\phi_{j}^{(t)}\sum\limits_{k\neq i}p_{ki}^{(t)}\sum\limits_{k\neq j}p_{kj}^{(t)}\Big]^{\frac{1}{2}}\Big\}}=\widetilde{\Gamma}_t,
\end{eqnarray}
where $\phi_{ij}^{(t)}=\max\{\phi_i^{(t)},\phi_j^{(t)}\}-\min\Big\{\frac{1}{a_{ii}-\sum\limits_{k\neq i}\frac{a_{ik}a_{ki}}{a_{kk}}},\frac{1}{a_{jj}-\sum\limits_{k\neq j}\frac{a_{jk}a_{kj}}{a_{kk}}}\Big\}.$
\end{thm}
\begin{proof}
Let $A^{-1}=[\alpha_{ij}]$. Since $A\in M_n$ is strictly diagonally dominant, we have, by Lemma \ref{zjx2015-1lemma1} and Lemma \ref{hzg}, that
\begin{eqnarray}\label{th5-equ2}
\frac{1}{a_{ii}-\sum\limits_{k\neq i}\frac{a_{ik}a_{ki}}{a_{kk}}}\leq \alpha_{ii}\leq \phi_i^{(t)}, i\in N.
\end{eqnarray}
Then
\begin{eqnarray}\label{th5-equ3}\nonumber
(\alpha_{ii}-\alpha_{jj})^2&\leq&\Bigg[\max\{\phi_i^{(t)},\phi_j^{(t)}\}-\min\bigg\{\frac{1}{a_{ii}-\sum\limits_{k\neq i}\frac{a_{ik}a_{ki}}{a_{kk}}},\frac{1}{a_{jj}-\sum\limits_{k\neq j}\frac{a_{jk}a_{kj}}{a_{kk}}}\bigg\}\Bigg]^2\\&=&[\phi_{ij}^{(t)}]^2.
\end{eqnarray}
By Theorem \ref{th2}, inequalities (\ref{th5-equ2}) and (\ref{th5-equ3}), we have
\begin{eqnarray*}
 \tau(A)&\geq& \frac{2}{\max\limits_{i\neq j}\Big\{\alpha_{ii}+\alpha_{jj}+\Big[(\alpha_{ii}-\alpha_{jj})^2+4\alpha_{ii}\alpha_{jj}\sum\limits_{k\neq i}p_{ki}^{(t)}\sum\limits_{k\neq j}p_{kj}^{(t)}
 \Big]^{\frac{1}{2}}\Big\}}\\
 &\geq& \frac{2}{\max\limits_{i\neq j}\Big\{\phi_i^{(t)}+\phi_j^{(t)}+\Big[(\phi_{ij}^{(t)})^2+4\phi_{i}^{(t)}\phi_{j}^{(t)}\sum\limits_{k\neq i}p_{ki}^{(t)}\sum\limits_{k\neq j}p_{kj}^{(t)}
 \Big]^{\frac{1}{2}}\Big\}}.
 \end{eqnarray*}
 The proof is completed.
\end{proof}

Similar to the proof of Theorem \ref{th5}, the following theorem is obtained easily.
\begin{thm}\label{th6}
If $A=[a_{ij}]\in M_n$ is strictly diagonally dominant, then for $t=1,2,\ldots$,
\begin{eqnarray}\label{th6-equ1}
 \tau(A)\geq \frac{1}{\max\limits_{i\in N}\Big\{\big(1+\sum\limits_{k\neq i}p_{ki}^{(t)}\big)\phi_i^{(t)}\Big\}}=\widetilde{\Omega}_t.
 \end{eqnarray}
\end{thm}

\begin{thm}\label{th7}
The sequence $\{\widetilde{\Gamma}_t\}$~($\{\widetilde{\Omega}_t\}$), $t=1,2,\ldots$ obtained
from Theorem \ref{th5} (Theorem \ref{th6}) is monotone increasing with an upper bound
$\tau(A)$ and, consequently, is convergent.
\end{thm}
\begin{proof}
By Lemma \ref{zjx2015-1lemma1}, we have $1>p^{(t)}_{ji}\geq p^{(t+1)}_{ji}\geq
0,j,i\in N,j \neq i, t=1,2,\ldots$. Then, by the definitons of $\phi_i^{(t)}$, it is easy to
see that the sequence $\{\phi_i^{(t)}\}$ is monotone decreasing.
Further, by the definition of $\phi_{ij}^{(t)}$, we know that the sequence $\{\phi_{ij}^{(t)}\}$ is also monotone decreasing.
Thus, $\{\widetilde{\Gamma}_t\}$ ($\{\widetilde{\Omega}_t\}$) is monotonically
increasing sequence. Hence, the
sequence $\{\widetilde{\Gamma}_t\}$ ($\{\widetilde{\Omega}_t\}$) is convergent.
\end{proof}

\begin{thm}\label{th8} Let $A=[a_{ij}]\in M_n$ with $a_{11}=a_{22}=\cdots=a_{nn}$, and $A^{-1}=[\alpha_{ij}]$ be doubly stochastic.
Then, for $t=1,2,\ldots$,
\begin{eqnarray*}
&&(a)~\Gamma_t\geq\Omega_t\geq\frac{1}{[1+(n-1)\rho(J_A)]\max\limits_{i\in N}{\alpha_{ii}}};\\
&&(b)~\Gamma_t\geq\frac{2}{\max\limits_{i\neq j}\big\{\alpha_{ii}+\alpha_{jj}+[(\alpha_{ii}-\alpha_{jj})^2+4(n-1)^2\alpha_{ii}\alpha_{jj}\rho(J_A)^2]^{\frac{1}{2}}\big\}}\geq \frac{1}{[1+(n-1)\rho(J_A)]\max\limits_{i\in N}{\alpha_{ii}}};\\
&&(c)~\widetilde{\Omega}_t\geq\frac{1}{[1+(n-1)d]\max\limits_{i\in N}{\varphi_{i}}};\\
&&(d)~\widetilde{\Gamma}_t\geq\frac{2}{\max\limits_{i\neq j}\big\{\varphi_{i}+\varphi_{j}+[\varphi_{ij}^2+4(n-1)^2\varphi_{i}\varphi_{j}d^2]^{\frac{1}{2}}\big\}}.
\end{eqnarray*}
\end{thm}
\begin{proof}
Since $A^{-1}$ is doubly stochastic, by Lemma \ref{ssj}, we have
$|a_{ii}|=\sum\limits_{j\neq i}|a_{ij}|+1=\sum\limits_{j\neq i}|a_{ji}|+1.$
Then for any $i\in N,$
$
r_i=\max\limits_{j\neq i}\bigg\{\frac{|a_{ji}|}{|a_{jj}|-\sum\limits_{k\neq j,i}|a_{jk}|}\bigg\}
   =\max\limits_{j\neq i}\big\{\frac{|a_{ji}|}{1+|a_{ji}|}\big\}
   =\frac{\max\limits_{j\neq i}|a_{ji}|}{1+\max\limits_{j\neq i}|a_{ji}|}.
$
Since $f(x)=\frac{x}{1+x}$ is an increasing function on $(0,+\infty)$, we have
$$
r_i=\frac{\max\limits_{j\neq i}|a_{ji}|}{1+\max\limits_{j\neq i}|a_{ji}|}
\leq\frac{\sum\limits_{j\neq i}|a_{ji}|}{1+\sum\limits_{j\neq i}|a_{ji}|}
=\frac{\sum\limits_{j\neq i}|a_{ji}|}{|a_{ii}|}=1-\frac{1}{|a_{ii}|},i\in N.
$$
Since
$
J_A=\left[\begin{array}{cccc}
0&-\frac{a_{12}}{a_{11}}&\cdots&-\frac{a_{1n}}{a_{11}}\\
-\frac{a_{21}}{a_{22}}&0&\cdots&-\frac{a_{2n}}{a_{22}}\\
\vdots&\vdots&\ddots&\vdots\\
-\frac{a_{n1}}{a_{nn}}&-\frac{a_{n2}}{a_{nn}}&\cdots&0\\
\end{array}\right]\geq 0,
$
then the $i$th row sum is
$d_i=\frac{\sum\limits_{j\neq i}|a_{ij}|}{|a_{ii}|}=1-\frac{1}{|a_{ii}|},i\in N.$
Further, from $a_{11}=a_{22}=\cdots=a_{nn}$, we have $d_i=d_j,i,j\in N,i\neq j.$
Hence, $\rho(J_A)=d=1-\frac{1}{|a_{ii}|}, i\in N.$
Combining with Lemma \ref{zjx2015-1lemma1}, we have that
$1>\rho(J_A)=d\geq r_i\geq p^{(t)}_{ji}\geq 0, i,j\in N,j\neq i,t=1,2,\ldots.$
Obviously,
\begin{eqnarray}\label{bjequ}
(n-1)\rho(J_A)=(n-1)d\geq \sum\limits_{j\neq i}p_{ji}^{(t)},~ \varphi_i\geq \phi_i^{(t)}, i,j\in N, i\neq j,t=1,2,\ldots.
\end{eqnarray}

From inequality (\ref{bjequ}) and Remark \ref{rmk1}, clearly, the conclusion (a) follows.
From inequality (\ref{bjequ}), Theorem 4.2 in \cite{lcq} and Theorem \ref{th2}, the conclusion (b) follows.
From inequality (\ref{bjequ}) and Theorem \ref{th6}, the conclusion (c) follows.

Since $\frac{1}{a_{ii}-\sum\limits_{k\neq i}\frac{a_{ik}a_{ki}}{a_{kk}}}\geq \frac{1}{a_{ii}},~\varphi_i\geq\phi_i^{(t)}, i\in N, t=1,2,\ldots$
then by the definitions of $\varphi_{ij}$ and $\phi_{ij}^{(t)}$, we have $\varphi_{ij}\geq \phi_{ij}^{(t)}, i,j\in N, i\neq j, t=1,2,\ldots.$
Further, from inequality (\ref{bjequ}) and Theorem \ref{th5}, the conclusion (d) follows.
\end{proof}

\section{Numerical examples}
In this section, several numerical examples are given to verify the
theoretical results.

\begin{example}{\rm Let
\[A =  \left[\begin{array}{rrrrrrrrrr}
27&-2&-4&-1&-3&-3&-4&-5&-1&-3\\
-2&34&-13&-2&-4&-2&-5&0&-3&-2\\
-3&-5&34&-6&-4&-3&-5&-2&-3&-2\\
0&-3&-4&38&-13&-4&-1&-4&-3&-5\\
-3&-3&-1&-11&41&-9&-2&-3&-4&-4\\
-3&-5&-2&-3&-6&35&-1&-5&-5&-4\\
-5&-2&0&-5&0&-7&34&-8&-1&-5\\
-1&-4&-3&-2&-5&-1&-9&32&-1&-5\\
-4&-4&-2&-4&-4&-3&-2&-1&33&-8\\
-5&-5&-4&-3&-1&-2&-4&-3&-11&37.9\\
\end{array} \right].\]
It is easy to verify that $A\in M_{10}$. Since $a_{10,10}=37.9<38=\sum\limits_{j\neq 10}|a_{10,j}|$, $A$ is not strictly
diagonally dominant and weakly chained
diagonally dominant. Hence inequalities (\ref{shi-equ}), (\ref{tgx-equ2}), (\ref{lcq-equ2}), (\ref{zjxcor1-equ1}), (\ref{th5-equ1}) and (\ref{th6-equ1}) can not be used to estimate the lower
bounds of $\tau(A)$. Numerical results obtained from Theorem 3.1 of \cite{tian}, Theorem 4.1 of \cite{lcq}, Theorem 4 of \cite{wf}, Theorem 3 of \cite{zjx2016mineigenvaule} and Theorem \ref{th2}, i.e., inequalities (\ref{tgx-equ1}), (\ref{lcq-equ1}), (\ref{wf-equ1}), (\ref{zjx2016-equ1}) and (\ref{th2-equ1}) are given in Table 1 for the total number of
iterations $T=10$. In fact, $\tau(A)=0.8873$.}
\end{example}

\begin{table}[bthp]\label{Table-1}
\caption{The lower upper of $\tau(A)$} \vspace{-1pt}
\begin{center}
\def\temptablewidth{0.7\textwidth}
{\rule{\temptablewidth}{0.05pt}}
\begin{tabular*}{\temptablewidth}{@{\extracolsep{\fill}}llllll}
Method &$t$&$\Upsilon_t$&Method&$t$&$\Gamma_t$\\\hline
Theorem 3.1 of \cite{tian}&&0.7195&&&\\
Theorem 4 of \cite{wf}&&0.7223&&&\\
Theorem 4.1 of \cite{lcq}&&0.7260&&&\\
Theorem 3 of \cite{zjx2016mineigenvaule}
&$t=1$&0.7380&Theorem 2&$t=1$&0.7905\\
&$t=2$&0.7870&&$t=2$&0.8328\\
&$t=3$&0.8123&&$t=3$&0.8569\\
&$t=4$&0.8231&&$t=4$&0.8659\\
&$t=5$&0.8289&&$t=5$&0.8708\\
&$t=6$&0.8319&&$t=6$&0.8737\\
&$t=7$&0.8336&&$t=7$&0.8749\\
&$t=8$&0.8344&&$t=8$&0.8754\\
&$t=9$&0.8349&&$t=9$&0.8757\\
&$t=10$&0.8351&&$t=10$&0.8759\\
\end{tabular*}
{\rule{\temptablewidth}{0.05pt}}
\end{center}
\end{table}

\begin{example}{\rm Let
\[A =  \left[\begin{array}{rrrrrrrrrr}
41&-12&-1&-5&-3&-3&-4&-4&-3&-3\\
-9&42&-15&-2&0&-4&0&-3&-4&-4\\
-1&-5&43&-13&-3&-3&-5&-4&-4&-4\\
-3&-5&-6&36&-9&-4&-3&-1&0&-4\\
-4&-3&-5&-2&34&-10&-2&-1&-4&-2\\
-3&-1&-4&-2&-1&37&-15&-5&-2&-3\\
-5&-2&-2&-2&-4&-2&35&-8&-5&-4\\
-5&-5&-1&-4&-5&-3&0&33&-6&-3\\
-5&-3&-4&-3&-3&-2&-2&-3&37&-11\\
-3&-5&-4&-2&-5&-5&-3&-3&-8&38.1\\
\end{array} \right].\]
It is easy to see that $A\in M_{10}$ is strictly diagonally dominant.
Next, we use only the entries of $A$ to give the lower bounds of $\tau(A)$.
Numerical results obtained from Theorem 4.1 of \cite{Shi}, Corollary 3.4 of \cite{tian}, Corollary 4.4 of \cite{lcq}, Corollary 1 of \cite{zjx2016mineigenvaule}, Theorem 14 of \cite{lcq2014}, and Theorem \ref{th5}, i.e., inequalities (\ref{shi-equ}), (\ref{tgx-equ2}), (\ref{lcq-equ2}), (\ref{zjxcor1-equ1}) and (\ref{th5-equ1}) are given in Table 2 for the total number of
iterations $T=10$. In fact, $\tau(A)=1.0987$.}
\end{example}

\begin{table}[bthp]\label{Table-2}
\caption{The lower upper of $\tau(A)$} \vspace{-1pt}
\begin{center}
\def\temptablewidth{0.7\textwidth}
{\rule{\temptablewidth}{0.05pt}}
\begin{tabular*}{\temptablewidth}{@{\extracolsep{\fill}}llllll}
Method &$t$&$\widetilde{\Upsilon}_t$&Method &$t$&$\widetilde{\Gamma}_t$\\\hline
Theorem 4.1 of \cite{Shi}&& 0.1000&&&\\
Corollary 3.4 of \cite{tian}&&0.1265&&&\\
Theorem 14 of \cite{lcq2014}&&0.1300&&&\\
Corollary 4.4 of \cite{lcq}&&0.1559&&&\\
Corollary 1 of \cite{zjx2016mineigenvaule}
&$t=1$&0.6219&Theorem \ref{th5}&$t=1$&0.6288\\
&$t=2$&0.8035&&$t=2$&0.8192\\
&$t=3$&0.9018&&$t=3$&0.9302\\
&$t=4$&0.9565&&$t=4$&0.9968\\
&$t=5$&0.9838&&$t=5$&1.0337\\
&$t=6$&0.9994&&$t=6$&1.0533\\
&$t=7$&1.0085&&$t=7$&1.0649\\
&$t=8$&1.0125&&$t=8$&1.0718\\
&$t=9$&1.0142&&$t=9$&1.0760\\
&$t=10$&1.0147&&$t=10$&1.0785\\
\end{tabular*}
{\rule{\temptablewidth}{0.05pt}}
\end{center}
\end{table}

\begin{rmk} Numerical results in Table 1 and Table 2 show that~:
\par (a) Lower bounds obtained from Theorem \ref{th2} and Theorem \ref{th5} are
bigger than these corresponding bounds in [3-8].

 (b)
These sequences obtained from Theorem \ref{th2} and Theorem \ref{th5} are monotone increasing.

(c) These sequences obtained from Theorem \ref{th2} and Theorem \ref{th5} approximates
effectively to the true value of $\tau(A)$.
\end{rmk}

\begin{example}\label{eg3}
\emph{Let $A=[a_{ij}]\in \mathbb{R}^{10\times 10}$, where $a_{ii}=10, a_{ij}=-1,i,j\in N,i\neq j.$
It is easy to see that $A\in M_{10}$ is strictly diagonally dominant. By Theorem \ref{th2}, Theorem \ref{th3} and Theorem \ref{th6} for $T=10$, respectively, we all have $\tau(A)\geq 1$ when $t=1$. In fact,
$\tau(A)=1$.}
\end{example}

\begin{rmk} Numerical results in Example \ref{eg3} show that the lower bounds obtained from Theorem \ref{th2}, Theorem \ref{th3} and Theorem \ref{th6} could reach the true value of $\tau(A)$ in some cases.
\end{rmk}

\section{Further work}
In this paper, we present several convergent sequences to approximate $\tau(A)$. Then an interesting problem is how accurately these bounds can be computed.
At present, it is very difficult for the authors to give the error analysis. We will continue to
study this problem in the future.

\section*{Competing interests}
The authors declare that they have no competing interests.

\section*{Author's contributions}
All authors contributed equally to this work. All authors read and
approved the final manuscript.

\section*{Acknowledgements}
\ifthenelse{\boolean{publ}}{\small}{}
This work is supported by the National Natural Science Foundation of
China (Nos.11361074,11501141), Foundation of Guizhou Science and
Technology Department (Grant No.[2015]2073).





\end{bmcformat}
\end{document}